\documentclass[11pt]{article}

\usepackage{amsmath}
\usepackage{amssymb}
\usepackage{eucal}

\begin{document}

\baselineskip 16pt

\title{On a  question of L.A. Shemetkov concerning 
   the   intersection\\ of   $\cal F$-maximal subgroups of finite
 groups\\-{\small Dedicated to 
 Professor K.P. Shum   on the occasion of his 70-th birthday}}

\author{ \\
{ Alexander  N. Skiba}\\
{\small Department of Mathematics,  Francisk Skorina Gomel State University,}\\
{\small Gomel 246019, Belarus}\\
{\small E-mail: alexander.skiba49@gmail.com}}

\date{}
\maketitle

\begin{abstract} We investigate the 
influence of the intersection of  the $\cal F$-maximal subgroups 
 on the structure of a finite group. In particular,
 answering   a question of  L.A Shemetkov we give conditions
 under which a hereditary saturated formation $\cal F$ has a property
 that for any finite  group $G$, 
 the ${\cal F}$-hypercentre of $G$ coincides with the intersection
of all   $\cal F$-maximal subgroups of $G$.

\end{abstract}

\let\thefootnoteorig\thefootnote

\renewcommand{\thefootnote}{\empty}


                                                                                                                                                                                                                                                                                                                                                                                                                                                                     \let\thefootnote\thefootnoteorig

\section{Introduction}

 Throughout this paper, all groups are finite.   We use 
$\cal N$, $\cal S$ and  $\cal U$  to denote the classes of all nilpotent,  of 
all soluble    and of all   supersoluble groups respectively.

Let $\cal X$ be a  class  of groups.  The symbol  
 $\pi ({\cal X})$ denotes the set of all primes $p$ such that 
$p$ divides $|G|$ for some $G\in {\cal X}$.    
A chief factor $H/K$ of a group
$G$ is called ${\cal X}$-central in $G$ provided $(H/K)\rtimes (G/C_{G}(H/K))\in
\cal X$ (see \cite[p. 127-128]{Shem-Sk}).  A normal subgroup $N$ of $G$ is 
said to be \emph{$\cal X$-hypercentral in $G$} if either $N=1$ or $N\ne 1$ and 
  every  chief factor of $G$  below $N$ is  ${\cal X}$-central in $G$. The symbol 
 $Z_{\cal X}(G)$   denotes the ${\cal X}$-hypercentre of $G$, that is, 
the product of all normal  ${\cal X}$-hypercentral subgroups of $G$ 
    \cite[p.~389]{DH}. If $1\in {\cal X}$ and $G$ is a group, then we write 
$G^{\cal X}$ to denote the intersection of all normal subgroups $N$ of $G$ 
with  $G/N\in {\cal X}$.   A group $G$ is called \emph{$s$-critical 
for ${\cal X}$} or simply \emph{${\cal X}$-critical } if $G$ is not in ${\cal X}$ but all proper subgroups of $G$ 
are in ${\cal X}$ \cite[p. 517]{DH}.   A subgroup $U$ of 
a group $G$ is called \emph{$\cal X$-maximal} in  $G$ provided that 
 (a)  $U\in {\cal X}$,   and   (b) if $U\leq V\leq G$  and $V\in {\cal X}$,
 then $U=V$ \cite[p. 288]{DH}.

Some  classes of  $\cal X$-maximal subgroups ($\cal X$-projectors,
 $\cal X$-injectors, $\cal X$-covering subgroups and  at al)  have 
been  studied by a large   number of authors and  they play an important role
 in  the theory of   soluble groups
 \cite{DH}.    In this paper, we investigate the 
influence of the intersection of all   $\cal X$-maximal subgroups of a 
group $G$  on the structure of  $G$.   We denote this intersection by 
  $\text{Int}_{\cal  X}(G)$.

In the paper \cite{BaerI}, Baer proved  that  
$\text{Int}_{\cal  N}(G)$ coincides with the  hypercentre
 $Z_{\infty }(G)=Z_{\cal N}(G)$ of  $G$. But in general,   
$Z_{\cal X}(G) < \text{Int}_{\cal  X}(G) $ even when  ${\cal X}={\cal U}$
 and $G$ is soluble (see  Example  5.13 below). 

L.A. Shemetkov asked in 1995 at the Gomel Algebraic Seminar the following 
question (the formulation of this question  was also given in \cite[p. 
41]{Sidorov}):  {\sl What  are the non-empty hereditary saturated formations
 $\cal F$ with the  property   that  for each group $G$, the equality 
 $$ \text{Int}_{\cal F}(G)=Z_{\cal F}(G) \eqno(*) $$  holds?}  
Our  main goal here is to give an  answer to this question.

A class   ${\cal F}$  of groups is said to be a \emph{formation} if either
 ${\cal F}=\varnothing$ or ${\cal F}\ne \varnothing$ and for any group $G$, each
 homomorphic  image of $G/G^{\cal F}$ belongs to ${\cal F}$. A formation ${\cal F}$
is said to be: \emph{saturated} if $G\in {\cal F}$ whenever $G/\Phi (G) \in {\cal 
F}$; \emph{hereditary} if $H\in {\cal F}$ whenever $H\leq G \in {\cal 
F}$.

Let  ${\cal F}$  be a saturated formation with  $\pi ({\cal F})\ne \varnothing$. 
Then for any $p\in \pi ({\cal F})$ we write ${\cal F}(p)$
 to denote  the  intersection of all formations
containing the set $\{ G/O_{p',  p}(G) \mid G\in {\cal F}\}$, and let  
 $F(p)$  denote   the  class of all groups  $G$ such that $ G ^{{\cal F}(p)}$
 is a $p$-group.

{\bf Remark.} We will show (see  Lemma 2.1 below) that  the function $f$ of 
the form   $$ 
f:\Bbb{P}\to \{ \text{group formations}\},$$ where $f(p)=F(p)$ for all 
$p\in \pi ({\cal F})$, and $f(p)=\varnothing $ for all $p\not \in \pi ({\cal 
F})$, is   the canonical local definition of ${\cal F}$   (see p. 361 in \cite{DH}).  
  Therefore, our notation $F(p)$ follows    the terminology  of  \cite[Chapter IV]{DH}.

{\bf Definition.}    Let
 ${\cal F}$  be a  hereditary saturated formation
 with   $\pi ({\cal F})\ne \varnothing$.   We say that  $\cal F$ satisfies:

(1)  The \emph{boundary condition} 
if for any $p\in \pi ({\cal F})$, $G\in {\cal F}$ whenever $G$ is an
 $F(p)$-critical group.

(2)   The \emph{boundary condition in the class of all  soluble groups} 
if for any $p\in \pi ({\cal F})$, $G\in {\cal F}$ whenever $G$ is a 
soluble  $F(p)$-critical group.

If ${\cal F}$ is a non-empty  formation with $\pi ({\cal F})= 
\varnothing$, then ${\cal F}=(1)$ is the class of all  groups $G$ with 
$|G|=1$,  and therefore 
for any group $G$ we have $Z_{\cal F}(G)= 1=\text{Int}_{\cal 
F}(G)$. For the general case, we prove 

{ \bf Theorem A. } {\sl   Let
 ${\cal F}$  be a  hereditary saturated formation
 with   $\pi ({\cal F})\ne \varnothing$. Equality  (*)  holds for each 
group $G$ if and only if $\cal F$  satisfies
 the boundary condition. }

{ \bf Theorem B. } {\sl    Let
 ${\cal F}$  be a  hereditary saturated formation
 with   $\pi ({\cal F})\ne \varnothing$. Equality  (*)  holds for each 
soluble group $G$ if and only if $\cal F$  satisfies
 the boundary condition    in the class  of all soluble   groups. }

The proofs of Theorems A and B rely on  
the following general facts on the subgroup  $\text{Int}_{\cal F}(G)$.

{\bf Theorem C. } {\sl Let
 ${\cal F}$  be a  non-empty hereditary saturated formation. Let $H$, $E$ be 
  subgroups of a group   $G$,    $N$  a normal subgroup of $G$ 
 and $I=\text{Int}_{\cal F}(G)$.
 }

(a) {\sl  $\text{Int}_{\cal F}(H)N/N\leq \text{Int}_{\cal F}(HN/N)$}.

(b) {\sl   $ \text{Int}_{\cal F}(H) \cap
E\leq  \text{Int}_{\cal F}(H\cap E)$.  }

(c) {\sl   If $H/H\cap I\in {\cal F}$, then $H\in {\cal F}$.}

(d)  {\sl If $H\in {\cal F}$, then $IH\in {\cal F}$}

(e) {\sl If $N\leq I$, then 
 $I/N= \text{Int}_{\cal F}(G/N)$}.

(f) {\sl  $\text{Int}_{\cal F}(G/I)= 1$}.

 (g) {\sl  If  every     $\cal F$-critical subgroup of $G$ 
 is soluble and  $ \psi _{e}(N) \leq
 I$, then $N\leq I$. }

(h)   $Z_{\cal F}(G) \leq I$.

It this theorem  $ \psi _{e}(N)$ denotes the subgroup of $N$ generated by all its 
cyclic subgroups of prime order and  order 4 \cite{La}.
 
We prove Theorems A, B and C  in Section 3.
In Section 4 it is shown that   
 the formation of all  nilpotent groups, the formation of all $p$-decomposable 
groups (for any prime $p$),
 and the formation of all   groups $G$ with  $G'\leq F(G)$  satisfy 
 the boundary condition,
 and that the formation of all soluble groups  of 
 nilpotent length  at most $r$ (for any fixed $r\in {\mathbb N}$) satisfies
 the boundary condition in the class of all soluble 
groups.
 We also  consider  here some classes
 of saturated   formations
 which do not satisfy  the boundary    condition.    
 Finally, in Section 5, some further applications 
 of the subgroup $\text{Int}_{\cal F}(G)$ are discussed.

 All unexplained notation and terminology are standard. The
reader is referred to   \cite{DH},  \cite{Guo}  and \cite{Bal-Ez}  if necessary.

\section{Preliminaries}

 The 
product ${\cal M}{\cal H}$ of the formations    ${\cal M}$  and  ${\cal H}$  is the 
class of all groups $G$ such that $G^ {\cal H}\in {\cal M}$. 
We use 
${\cal G}_{\pi}$  to denote the class of all $\pi$-groups. In particular,  
we write ${\cal G}_{p}$  to denote the class of all $p$-groups if ${\pi}=
\{p\}$, $p$ is a prime.  The product of any two formations  is itself a
 formation \cite[Chapter IV, Theorem 1.8]{DH}.  Therefore, if $\cal F$ is 
a saturated formation and if $p\in \pi ({\cal F})$, then $F(p)= {{\cal 
G}_{p}}{\cal F}(p)$ is a formation.

A function $     f:\Bbb{P}\to \{ \text{group formations}\}$ is called a
  \emph{formation function}.  The symbol $LF(f)$ denotes the collection of all   
groups $G$ such that either $G=1$ or $G\ne 1$ and $G/C_{G}(H/K)\in f(p)$  
for every chief factor $H/K$ of $G$ and every $p\in \pi (H/K)$. A 
formation function $f$ is called  \emph{integrated} if $f(p)\subseteq LF(f) $
 for all
primes $p$, and  \emph{full} if $f(p)={\cal G}_{p}f(p)$  for all primes
 $p$.    If for a 
formation $\cal F$ we have  ${\cal F}=LF(f)$, then $f$ is called a local 
definition of $\cal F$.   
 It is well known that $O_{p', p}(G)= \cap \{C_{G}(H/K) \mid H/K
  \  \text{is a chief factor of }\     G   \  \text{and} \  p \in \pi (H/K) \}$.
 Therefore,   
 $G\in {\cal F}=LF(f)$ if and only if either $G=1$ or $G\ne 1$ and
 $G/O_{p', p}(G)\in f(p)$ for all $p\in \pi (G)$.

 { \bf Lemma 2.1.}  {\sl  Let $\cal F$ be a non-empty saturated formation. Then  
${\cal F}=LF(f)$, where $f(p)=F(p)\subseteq {\cal F}$ for all 
 $p\in \pi ({\cal F})$, and $f(p)=\varnothing $ for all primes $p\not \in
 \pi ({\cal F})$. }

{\bf Proof.}  Define a function $t$ 
as follows:

$$t(p)= \begin{cases} {\cal F}(p), \ &\text{if} \ p\in \pi ({\cal F}),
 \\  \varnothing , &\text{if} \ p\not \in \pi ({\cal F})

\end{cases}.$$ Let ${\cal M}=LF(t)$. Then ${\cal F}\subseteq {\cal M}$. On 
the other hand, by the  Gasch\"utz-Lubeseder-Schmid theorem 
  \cite[Chapter IV, Theorem 4.4]{DH}, there is a 
 formation function $h$ such that ${\cal F}=LF(h)$.  Moreover, 
 $t(p)\leq h(p)$ for all primes $p$ and therefore ${\cal 
M}  \subseteq {\cal F}$. Hence ${\cal F}={\cal M}=LF(t)$. Now the assertion 
 follows from Proposition 3.8 (a) in \cite[Chapter IV]{DH}.

From   Theorem 17.14 in \cite{Shem-Sk} we  
get

 { \bf Lemma 2.2. } {\sl  Let $\cal F$ be a  non-empty saturated formation.
   A chief factor $H/K$ of a group $G$ 
  is ${\cal F}$-central in $G$ if and only if $G/C_{G}(H/K)\in F(p)$
 for all primes  $p\in \pi (H/K)    $.}

In view of   Lemma 2.1  and Proposition 3.16     in \cite[IV]{DH}
  we have  

 { \bf Lemma 2.3. } {\sl  Let
 ${\cal F}$  be a  hereditary saturated formation. 
   Then for any prime $p\in \pi ({\cal F})$,
 the formation $F(p)$ is hereditary.}

We shall  need in our proofs a few facts about  the   ${\cal 
F}$-hypercentre.

 { \bf Lemma 2.4. } {\sl  Let
 ${\cal F}$  be a  non-empty  saturated formation. 
 Let $G$ be a group and $H\leq G$.}

(1) {\sl If $H$ is normal in $G$, then $Z_{\cal F}(G)H/H\leq Z_{\cal F}(G/H)$}

(2) {\sl If ${\cal F}$  is hereditary, then $Z_{\cal F}(G)\cap H\leq Z_{\cal F}(H)$.}

(3) {\sl If $G/Z_{\cal F}(G)\in {\cal F}$, then $G\in {\cal F}$.}

{\bf Proof.}  (1) This follows from the  $G$-isomorphism $Z_{\cal F}(G)H/H\simeq
 Z_{\cal F}(G)/Z_{\cal F}(G)\cap H$ since for any two $G$-isomorphic chief factors
 $H/K$ and   $T/L$ of $G$ we have  $(H/K)\rtimes (G/C_{G}(H/K))\simeq
 (T/L)\rtimes (G/C_{G}(T/L))$.

(2) Let $1= Z_{0} < Z_{1} < \ldots <  Z_{t} = Z_{\cal F}(G)$ be a chief
 series of $G$      below $Z_{\cal F}(G)$ and  $C_{i}= 
C_{G}(Z_{i}/Z_{i-1})$. Let $p$ be a prime divisor
 of $|Z_{i}\cap H/Z_{i-1}\cap H|=|Z_{i-1}(Z_{i}\cap H)/Z_{i-1}|$.  Then 
$p$ divides $|Z_{i}/Z_{i-1}|$, so $G/C_{i}\leq F(p)$ by Lemma 2.2. Hence 
by Lemma 2.3, $H/H\cap C_{i}\simeq C_{i}H/C_{i}\in F(p)$. But  $H\cap C_{i}\leq
 C_{H}(Z_{i}\cap H/Z_{i-1}\cap H)$. Hence
 $H/C_{H}(Z_{i}\cap H/Z_{i-1}\cap H)\in F(p)$ for all primes $p$ dividing 
$|Z_{i}\cap H/Z_{i-1}\cap H|$. Thus  $Z_{\cal F}(G)\cap H\leq Z_{\cal 
F}(H)$   by Lemma  2.2 and \cite[Chapter A, Theorem 3.2]{DH}.

(3) This follows from Lemmas 2.1, 2.2 and the Jordan-H\"older theorem 
\cite[Chapter A, Theorem 3.2 ]{DH}.

The following lemma is a corollary of general  results on  $f$-hypercentral action
 (see \cite[Chapter IV, Section 6]{DH}). For reader's 
convenience, we give a direct proof.

{\bf Lemma 2.5.} {\sl  Let
 ${\cal F}$  be a   saturated formation.  Let $E$ be a normal
 $p$-subgroup of a group $G$. 
 If  $ E\leq Z_{\cal F}(G)$, then
 $G/C_{G}(E)\in F(p)$.}

{\bf Proof.}   Let  $1=E_{0} < E_{1} <  \ldots  < E_{t} =E$   be a chief series of
 $G$ below $E$. Let  $C_{i} 
=C_{G}(E_{i}/E_{i-1})$ and    $C= C_{1}\cap \ldots \cap C_{t}$. 
Then  $C_{G}(E)\leq C$ and   so 
 $C/C_{G}(E)$ is a $p$-group by Corollary 3.3 in \cite[Chapter 5]{Gor}. On the other
 hand, by Lemma 2.2, 
   $G/C_{i}\in F(p)$, so $G/C\in F(p)$. Hence  
 $G/C_{G}(E)\in F(p)={\cal G}_{p}F(p)$.

{\bf Lemma 2.6.} {\sl  Let  $G$ be a group and $p$ a prime such that $O_{p}(G)=1$.
If  $G$ has a unique minimal normal subgroup, then there exists  a  simple 
   ${\mathbb F}_{p}G$-module
 which is     faithful  for $G$.  }

{\bf Proof.}  Let $C_{p}$ be a group of order $p$. Consider 
 $A=C_{p}\wr G=K\rtimes G$,  the regular wreath product of $C_{p}$ with $G$,
  where  $K$ is the base group of  $A$.
 Let $$1=K_{0} < K_{1}  < 
\ldots <  K_{t}=K,
\eqno(*)$$

where  $K_{i}/K_{i-1}$ is a chief  
factor of $A$ for all  $i=1,  \ldots , t$. Let 
 $C_{i}=C_{A}(K_{i}/K_{i-1})$,  $N$  a minimal normal subgroup of $G$
 and $C= C_{1}\cap  \ldots \cap
  C_{t} $.
 Suppose that $C_{i}\cap G\ne 1$ for 
all $i=1,  \ldots , t$. Then $N\leq C\cap G$. Hence $N$ stabilizes Series (*),
 so $N$ is a $p$-group by Corollary 3.3 in \cite[Chapter 5]{Gor},
 which implies $N\leq O_{p}(G)$.    This contradiction shows that
 for some $i$ we have  
$C_{G}(K_{i}/K_{i-1})=1$.   The lemma is proved.

{\bf Lemma 2.7.}   Let
 ${\cal F}$  be a  non-empty  saturated formation.

 (1) {\sl  If for some prime $p$ we have 
 $ {\cal F}={\cal G}_{p}{\cal F}$, then  $F(p)={\cal F}$.
}

(2) {\sl If  ${\cal F}={\cal N}{\cal H}$ for some non-empty
 formation ${\cal H}$, then 
 $F(p)= {\cal G}_{p}{\cal H}$ for all primes $p$. } 

{\bf Proof.} (1)  By Lemma  2.1, $ F(p)\subseteq {\cal F}$, so  we need only
  prove that $ {\cal 
F}\subseteq F(p)$. Suppose that this is false and let $A$ be a group of 
minimal order in ${\cal F}\diagdown F(p)$. Then $A^{F(p)}$ is a unique
 minimal normal subgroup of $A$  and $O_{p}(A)=1$. By Lemma 2.6  
  there is a  simple  
 ${\mathbb F}_{p}A$-module $P$    which is     faithful  for $A$.
Then $G=P\rtimes  A \in {\cal G}_{p}{\cal F}= {\cal F}$, so
 $A \simeq G/P=G/O_{p', p}(G) \in F(p)$, a contradiction. Thus   
$F(p)={\cal F}$.

(2)  The inclusion $F(p)\subseteq {\cal 
G}_{p}{\cal H}$ is evident.  Suppose that ${\cal G}_{p}{\cal 
H}\not \subseteq F(p)$ and let $A$ be a group of minimal order
 in ${\cal N}_{p}{\cal 
H}\diagdown F(p)$.   Then $A^{F(p)}$ is a unique minimal normal subgroup of 
$A$ and $O_{p}(A)=1$.  Hence  $A\in {\cal H}$ and  there exists a  simple  
 ${\mathbb F}_{p}A$-module $P$    which is     faithful  for $A $.
Then $G=P\rtimes  A \in {\cal G}_{p}{\cal H}\subseteq {\cal F}$, so
 $A \simeq G/P=G/O_{p', p}(G) \in F(p)$, a contradiction.    The lemma is proved.

 { \bf Lemma 2.8} \cite[Chapter VI, Theorem 25.4]{26}. 
 {\sl  Let
 ${\cal F}$  be a   saturated formation. Let   $G$ be a group whose
 $\cal F$-residual
  $G^{\cal F}$ is soluble. Suppose that every maximal subgroup of $G$ not
 containing $G^{\cal F}$ belongs to ${\cal F}$.}

(a) {\sl   $P=G^{\cal F}$ is a $p$-group for some prime $p$ and $P$  is of
 exponent $p$ or of exponent  $4$ (if $P$ is a
 non-abelian $2$-group)}.

(b) {\sl   $P/\Phi (P)$  is a chief factor of $G$  and   $(P/\Phi (P))
 \rtimes  (G/C_{G}(P/\Phi (P)))\not \in {\cal F}$.}

Let  $H$ and $K$ be  subgroups of a group $G$.  If $HK=G$, then $K$ is 
called  a \emph{supplement} of $H$ in $G$. If, in addition, $HT\ne G$ for all
 proper subgroups $T$ of $K$, then   $K$ is called a \emph{minimal supplement} 
of $H$ in $G$.

 { \bf Lemma 2.9. } {\sl  Let
 ${\cal F}$  be a  hereditary saturated formation.  Let $N\leq U\leq G$, where $N$ is a normal subgroup of a 
group $G$. }

(i) {\sl If $G/N\in {\cal F}$ and $V$ is a minimal supplement of $N$ in $G$, then 
 $V\in {\cal F}$.}

(ii) {\sl  If $U/N$ is an $\cal F$-maximal subgroup of $G/N$, then  $U=U_{0}N$ 
for some $\cal F$-maximal subgroup $U_{0}$ of $G$}.

(iii)  {\sl      If $V$ is an  $\cal F$-maximal 
  subgroup of $U$, then $V=H\cap U$ for 
some   $\cal F$-maximal   subgroup $H$ of $G$. }

 {\bf Proof.} (i)  It is clear that $V\cap N \leq \Phi (V)$. Hence from 
   $V/V\cap N \simeq VN/N=G/N \in {\cal F}$ we have  $V\in {\cal F}$ since
 ${\cal F}$ is  saturated.

 (ii) 
Let $V$ be a minimal supplement of $N$ in
 $U$.  Then  $V\in {\cal F}$ by (i). 
 Let $U_{0}$ be an $\cal F$-maximal subgroup of $G$  such that $V\leq  
U_{0}$. Then $U_{0}N/N\simeq U_{0}/U_{0}\cap N \in {\cal F}$  and 
  $U/N\leq U_{0}N/N$. Hence 
$U=U_{0}N$.

(iii)    Let  $H$ be an  $\cal F$-maximal   subgroup  of $G$ such that 
$V\leq H$. Then  $V\leq H\cap U\in {\cal F}$ since ${\cal F}$ is hereditary,
 which implies  $V=H\cap U$.

 {\bf Lemma 2.10. }  {\sl 
Let  ${\cal F}$  be a saturated formation with   
 $p\in \pi ({\cal F})$. Suppose that $G$ is a group of minimal order
 in the set of all  
  $F(p)$-critical  groups $G$ with $G\not \in {\cal F}$.
 Then  $O_{p}(G)=1=\Phi (G)$ and 
$G^{\cal F}$ is  a unique minimal normal    subgroup of $G$.
  }

{\bf Proof.} Let $N$ be a minimal normal subgroup of $G$. Then $G/N\in 
{\cal F}$. Indeed, suppose that    $G/N\not \in     {\cal F}$.
 Since $F(p)$ is a formation, and $F(p)\subseteq {\cal F}$  by Lemma 2.1, 
it follows that  $G/N \not \in F(p)  $ and that  every 
maximal subgroup of $G/N$   belongs to $F(p)$. Thus  $G/N$  is an 
  $F(p)$-critical  group with $G/N\not \in     {\cal F}$. But then $|G/N| < |G|$ 
 contradicts the minimality of $G$. Hence  $G/N \in     {\cal F}$.
Since ${\cal F}$ is a saturated formation,  $N=G^{\cal F}$ is a unique
 minimal normal   subgroup of $G$ and  $\Phi (G)=1$. Suppose that  $N\leq  O_{p}(G)$  and 
let $M$ be a maximal subgroup of $G$ such that $G=NM$. Then $G/N\simeq 
M/N\cap M \in F(p)={\cal G}_{p}F(p)$, so $G\leq F(p)\subseteq {\cal F}$. 
This contradiction shows that $O_{p}(G)=1$.

 {\bf Lemma 2.11 }\cite[Chapter 1, Lemma 4.4]{26}. 
 {\sl Let $L$ be a normal subgroup of a group $G$ such that $L\leq \Phi 
(G)$. If $G/L$ has a normal Hall $\pi $-subgroup, then does $G$.  

 }

{\bf Lemma 2.12} \cite[Lemma 2.2]{Sk}. {\sl Let $p\ne q$ be primes
dividing  the order of a group $G$, $P$ a Sylow $p$-subgroup of $G$. 
 If every maximal
subgroup of $P$ has a $q$-closed supplement in $G$, then $G$ is
$q$-closed.}

The following lemma is well known.

{\bf Lemma  2.13.} {\sl  Let $A$ and $B$ be  proper subgroups of $G$
such that $G=AB$. Then  $A^{x}B=G$  and  $G\ne AA^{x}$ for all
$x\in G$.}

\section{Proofs of Theorems A, B and  C}

{\bf Proof of Theorem C}         
(a)  First we suppose that $H=G$.
  If $U/  N$ is an $\cal F$-maximal subgroup of $G/N$, then for some 
 $\cal F$-maximal subgroup $U_{0}$ of $G$ we have 
   $U=U_{0}N$ by Lemma 2.9 (ii). Let $\text{Int}_{\cal F}(G/N)= U_{1}/N
  \cap  \ldots 
\cap U_{t}/N$, where  $U_{i}/N$ is an $\cal F$-maximal subgroup of $G/N$ 
for all $i=1,  \ldots , t$. Let $V_{i}$ be an $\cal F$-maximal subgroup 
of $G$ such that  $U_{i}=V_{i}N$. Then $I\leq  V_{1}  \cap \ldots 
\cap V_{t}$, so  
$IN/N\leq \text{Int}_{\cal F}(G/N)$.

 Now let $H$ be any subgroup of 
$G$. And let  $f:H/H\cap N \to HN/N$  be the canonical isomorphism from
$H/H\cap N $ onto $HN/N$. Then $f(\text{Int}_{\cal F}(H/H\cap N))=\text{Int}_{\cal
F}(HN/N)$ and $f(\text{Int}_{\cal F}(H)(H\cap N)/(H\cap N))=\text{Int}_{\cal
F}(H)N/N$. But  from   above  we have $\text{Int}_{\cal F}(H)(H\cap N)/(H\cap N)\leq
\text{Int}_{\cal F}(H/H\cap N)$. Hence $\text{Int}_{\cal F}(H)N/N\leq
 \text{Int}_{\cal
F}(HN/N)$.

(b) If $V$ is any $\cal F$-maximal  subgroup of $H$, then $V=H\cap U$ for 
some
 $\cal F$-maximal   subgroup $U$ of $G$ by Lemma 2.9 (iii).
 Thus there are  $\cal F$-maximal
 subgroups $U_{1},   \ldots , U_{t}$ of $G$ 
 such that  $\text{Int}_{\cal F}(H)=U_{1} \cap   \ldots \cap 
 U_{t}\cap H$,
hence    $I\cap H\leq  \text{Int}_{\cal F}(H)$ and   
 $ \text{Int}_{\cal F}(H) \cap
E     =  \text{Int}_{\cal F}(H) \cap
(H\cap E)\leq \text{Int}_{\cal F}(H\cap E)$. 

  (c) First we suppose that $H=G$. 
   Let $U$ be a minimal supplement of $I$ in $G$. Then 
 $U\in {\cal F}$ by Lemma 2.9 (i). 
  Let $V$ be an   $\cal F$-maximal subgroup of $G$ 
containing  $U$. Then  $G= IU= \leq V\in {\cal F}$. Finally, in
 the general case  
 we have  $I\cap H\leq  \text{Int}_{\cal F}(H)$ by  (b), so 
from  $H/H\cap I\in {\cal F}$ we
deduce  $H/ \text{Int}_{\cal F}(H)\in {\cal F}$ and hence 
$H  \in {\cal F}$.

(d)  Since $H\in {\cal F}$,  $HI/I \simeq H/H\cap I\in {\cal F}$. By (b), 
$I    \leq \text{Int}_{\cal F}(HI)$. Hence
 $HI/\text{Int}_{\cal F}(HI) \in {\cal F}$.
 Thus $HI\in {\cal F}$ by (c).

(e)       In view of Lemma 2.9 (ii) it is enough to prove that if  $U$ is an
 $\cal F$-maximal 
subgroup of $G$, then  $U/N$  is an $\cal F$-maximal
 subgroup of $G/N$.  Let  $U/N\leq X/N$, where $X/N$  is an
 $\cal F$-maximal subgroup of $G/N$. By Lemma 2.9 (ii),  $X=U_{0}N$ 
for some 
$\cal F$-maximal subgroup $U_{0}$ of $G$. But $N\leq U_{0}$, so  $U/N\leq 
U_{0}/N$ and hence $U=U_{0}$. Thus $U/N= X/N$.

(f) This follows from (e).

(g) Suppose that this assertion is false and let $G$  be  a counterexample 
with $|G||N|$ minimal. Then there is an $\cal F$-maximal  subgroup $U$  
of $G$ 
such that $N\nleq U$.  Let $E=NU$. Then $E/N\simeq U/U\cap N\in {\cal F}$.
  By (b), $\psi _{e}(N)\leq I\cap E  \leq  \text{Int} _{\cal F}(E)$.  
Suppose that $E\ne G$. Then  $ N\leq \text{Int} _{\cal F}(E)$  by the choice 
of $(G, N)$, so $G/\text{Int} _{\cal F}(E) \in {\cal F}$.  Hence $E\in {\cal F}$ by (c),
 so $U=E$. Therefore $N\leq U$, a 
contradiction. Thus $E=G$. Let $M$ be any maximal subgroup of $G$. We show  
that $M\in {\cal F}$.   Since 
 $ \psi _{e}(N\cap M)\leq \psi _{e}(N)$,   $ \psi _{e}(N\cap M)\leq 
 I\cap M$. Hence $ \psi _{e}(N\cap M)   \leq  \text{Int}_{\cal F}(M)$ 
by (b). Therefore  $N\cap M\leq  \text{Int}_{\cal F}(M)$ by the choice of $(G, N)$.
 Note also that $M/M\cap N\in
{\cal F}$. Indeed, if $N\leq M$, then $M/N\leq G/N\in {\cal 
F}$.  On the other hand, if $N\not \subseteq M$, then $M/M\cap N\simeq
 NM/N=G/N\in {\cal F}$ since ${\cal F}$ is hereditary. 
Therefore $M\in {\cal F}$ by (c).  Hence $I=\Phi (G)$  and $G$
 is an $ {\cal F}$-critical group. Since  
$G/N\in {\cal F}$, $\psi _{e}(G^{\cal F})\leq \psi _{e}(N)\leq
I$. Thus   for any 
$x\in G^{\cal F} \diagdown
 \Phi (G^{\cal F})$ we have $x\in \psi_{e}(N)\leq I=\Phi (G)$ by
 Lemma 2.8.
 Therefore 
$G^{\cal F}\leq \Phi (G)$, so $I=G\in {\cal F}$, a contradiction. 
Hence we have (g).

(h)   Let $H$  be  a subgroup of $G$ such that  $H\in { \cal F}$. Then
  $HZ_{\cal F}(G) /Z_{\cal F}(G) \simeq H/H\cap Z_{\cal F}(G) \in {\cal F}$ and 
 $Z_{\cal F}(G) \leq Z_{\cal F}(HZ_{\cal F}(G))$ by Lemma 2.4 (2).
 Hence $HZ_{\cal F}(G) \in {\cal F}$
 by Lemma 2.4  (3). Thus  $Z_{\cal F}(G) \leq  I$.

{\bf Proof of Theorem A.}     
First we  suppose that $\cal F$ satisfies  the boundary condition.
 We shall show  that  for  every group $G$ we have 
 $Z_{\cal F}(G)= \text{Int}_{\cal 
F}(G)$.  Suppose that this is 
false  and let $G$ be a counterexample with minimal order. Let $Z=Z_{\cal F}(G)$ 
and $I= \text{Int}_{\cal F}(G)$.
  Then    $Z < I$ by Theorem C (h),    so $I\ne 1$ and  $G\not 
\in {\cal F}$.  
 Let $N$ be a 
minimal normal subgroup of $G$, $L$  a minimal normal subgroup of $G$ contained
 in $ I$.

(1) {\sl   $IN/N\leq Z_{\cal F}(G/N)= \text{Int}_{\cal F}(G/N)$.}

Indeed, by Theorem C (a) we have $IN/N\leq  \text{Int}_{\cal F}(G/N)$. 
On the other hand, by the choice of $G$,  $\text{Int}_{\cal 
F}(G/N)=Z_{\cal F}(G/N)$.

(2) {\sl  $L\nleq Z$}.

Suppose that  $L\leq Z$. Then $Z/L=Z_{\cal F}(G/L)$ and $I/L=
 \text{Int}_{\cal F}(G/L)$ by
 Theorem C (e).  But by (1), $Z_{\cal F}(G/L)= \text{Int}_{\cal F}(G/L)$. 
Hence  $I /L=Z/L$, so $I=Z$, a contradiction.

(3) {\sl If $L\leq M < G$, then $L\leq Z_{\cal F}(M)$}.

Let $V$ be any  $\cal F$-maximal subgroup of $M$.  Then  $V=H\cap M$ for 
some   $\cal F$-maximal   subgroup $H$ of $G$ by Lemma 2.9 (iii).
 Hence $L$ is 
contained in the  intersection of all  $\cal F$-maximal subgroups  of 
$M$. But $|M| < |G|$, so $ \text{Int}_{\cal F}(M)=Z_{\cal F}(M)$  by the 
choice of $G$.   Hence  $L\leq Z_{\cal F}(M)$ 

(4) {\sl $L=N$ is a unique minimal normal subgroup of $G$}.

Suppose that $L\ne N$. From  Theorem C (a) and (1) we deduce that 
 $NL/N\leq Z_{\cal F}(G/N)$, 
so  from the $G$-isomorphism $NL/N\simeq L$ we obtain $L\leq Z$, which 
contradicts (2). 

(5) $L\nleq \Phi (G)$.

Suppose that $L\leq \Phi (G)$. Then $L$ is a $p$-group for some prime $p$. 
 Let $C=C_{G}(L)$. Let $M$ be any maximal 
subgroup of $G$. Then $L\leq M$, so  $L\leq Z_{\cal F}(M)$ by (3).
 Hence 
$M/M\cap C\in F(p)$ by Lemma 2.5. 
 If  $C\nleq M$, then  $G/C=CM/C\simeq M/M\cap C\in F(p)$, so 
 $L\leq 
Z_{\cal F}(G)$ by Lemma 2.2, contrary to (2). Hence $C\leq M$ for
 all maximal subgroups $M$
 of $G$, so $C$ is nilpotent. Therefore in view of (4), $C$ is a $p$-group
 since $C$       is normal
 in $G$. Hence  for every maximal subgroup $M$ of $G$
 we have $M\in {\cal G}_{p}F(p)=F(p)$. By Lemma 2.1,  $F(p)\subseteq {\cal 
F}$. Hence $G\not \subseteq F(p)$ and so $G$ is an  $F(p)$-critical group. 
 But $\cal F$ satisfies  the boundary condition and so  $G\in {\cal F}$, a contradiction.
Hence we have (5). 

(6) {\sl $L$ is not abelian}.

Suppose that $L$ is abelian. Then  from (4) and (5) we deduce that
  $G=L\rtimes M$ for some maximal subgroup $M$ of 
$G$ and  $C=C_{G}(L)=L.$
Let $E$ be a maximal subgroup of $M$, $V=LE$. Then  by (3), $L\leq Z_{\cal F}(V)$, 
so $E\simeq V/L=V/C_{V}(L)\in F(p)$ by Lemma 2.5.   
 Hence $M\in {\cal F}$ since $\cal F$ satisfies  the boundary condition. 
But $L\leq I$, so  $G\in {\cal F}$ by Theorem A (c), a contradiction.

{\sl The final contradiction for the sufficiency}.

Let  $p\in \pi (L)$. First we show that each  maximal  subgroup $M$  of $G$
 containing $L$ 
belongs to $F(p)$. By (3),  $L\leq 
Z_{\cal F}(M)$.  Let $$1=L_{0} < L_{1}  < 
\ldots <  L_{n}=L \eqno(*)$$ be a chief series of $M$ below $L$. Let
 $C_{i}=C_{M}(L_{i}/L_{i-1})$ and $C=C_{1}\cap \ldots \cap 
C_{n}$. Since by Lemma 2.2, $M/C_{i}\in F(p)$ for all $i=1,  \ldots , n $, $M/C\in F(p)$. 
By (4), $L$ is a unique minimal normal subgroup of $G$ and $L$ is 
non-abelian by (6). Hence $C_{G}(L)=1$, so for any minimal
 normal subgroup $R$ of $M$ we have
$R\leq L$. Suppose that $C\ne 1$ and let 
 $R$ be a minimal normal subgroup of $M$ contained in $C$. Then  $R\leq 
L$ and  $R\leq C_{A}(H/K)$ for each  chief factor $H/K$  of $M$ 
 by \cite[Chapter A, Theorem 3.2]{DH}. Thus 
 $R$ is abelian and hence $L$ is abelian. This contradiction shows that 
 $C=1$, so $M\in F(p)$.

 Now let $U$ be a minimal supplement of $L$ in $G$, $V$ a maximal 
subgroup of $U$. Then $LV\ne G$, so $LV\leq T$ for some maximal subgroup  $T$ of $G$.
Hence  $T\in F(p)$, so $V\in F(p)$ by Lemma 2.3. Therefore every maximal subgroup 
of $U$ belongs to $F(p)\subseteq {\cal F}$. Hence $U\in {\cal F}$, so $G\in {\cal F}$ by 
Theorem C (c).
This contradiction completes the proof of the  sufficiency.

Now  suppose that  the equality   $Z_{\cal F}(G)= \text{Int}_{\cal 
F}(G)$ holds   for    each group $G$.
  We shall show    that   $\cal F$  satisfies the boundary condition.
    Suppose that this is false.  Then there is  a prime $p\in \pi ({\cal 
F})$ such that the set of all $F(p)$-critical groups $A$ with $A\not \in 
{\cal F}$  is non-empty. Let us choose in this set a group $G$     with minimal 
$|G|$.
 Then  by Lemma 2.10, 
 $
G^{\cal F}$
 is a unique minimal normal subgroup of $G$ and  $O_{p}(G)=1=\Phi (G)$.
Hence by  Lemma 2.6, there exists  a    a  simple  
 ${\mathbb F}_{p}G$-module  $P$    which is     faithful  for $G$.   
 Let $A=P\rtimes G$ and  $M$ be
 any maximal subgroup 
of $A$. If $P\nleq M$, $M\simeq A/P\simeq G \not \in {\cal F}$. On the 
other hand, if $P\leq M$, $M=M\cap PG=P(M\cap G)$, where $M\cap G$ is a 
maximal subgroup of $G$. Hence   $M\cap G\in F(p)$, so
 $M\in {\cal G}_{p}F(p)=F(p)\subseteq   {\cal F}$ by Lemma 2.1.   
Therefore $P$ is contained in the intersection of all $\cal F$-maximal 
subgroups of $A$. Hence $P\leq Z_{\cal F}(A)$ by our assumption about $\cal F$, so 
$G\simeq A/P=A/C_{A}(P)\in F(p)\subseteq {\cal F}$  by Lemma 2.5. This 
contradiction completes the proof of the result.

{\bf Proof of Theorem B.}  See the proof of Theorem A.

\section{Some classes of formations satisfying the boundary condition }

{\bf Classes of soluble  groups with limited  nilpotent  length.}  
Following 
\cite[Chapter VII, Definitions 6.9]{DH} we write  $l(G)$ to denote the 
nilpotent length of the group $G$.  Recall that
 ${\cal N}^{r}$ is the product of $r$ copies of ${\cal N}$; ${\cal N}^{0}$
 is the class of groups of order 1 by definition. It is well known that
 ${\cal N}^{r}$
is the class of all soluble groups $G$ with  $l(G) \leq r$. It is known 
also  that  ${\cal N}^{r}$ is a hereditary saturated formation (see, for example,  
\cite[p. 358]{DH}).

{\bf Proposition 4.1.} {\sl For any $r\in {\mathbb N}$, the formation  ${\cal 
N}^{r}$   satisfies the boundary condition in the class of all soluble groups.
 The formation ${\cal N}$ satisfies the boundary condition. }

{\bf Proof.}  We proceed by induction on $r$. 
 Let ${\cal F}={\cal N}^{r}$, ${\cal H}={\cal N}^{r-1}$.  It is clear  
 that ${\cal F}={\cal N}{\cal H}$, so $F(p)={\cal N}_{p}{\cal H}$ for all primes
 $p$ by  Lemma 2.7 (2).  If $r=1$, 
then for any prime $p$ we have $F(p)={\cal N}_{p}$, so ${\cal F}={\cal N}$
 satisfies the boundary condition.

 Now suppose that $r > 1$. 
 Assume that  ${\cal F}$ does not   satisfy  
the boundary condition  in the class of all soluble groups. Then there is a  
prime $p$ such that  the set of all soluble $F(p)$-critical   groups $A$ with
 $A\not \in    {\cal F}$    is non-empty. Let $G $ be a group of minimal
 order in this set.  Then
  $O_{p}(G)=1=\Phi (G)$ and   $R=G^{\cal F}$ is  a unique minimal normal
 subgroup of $G$ by Lemma 2.10.  Hence  $G$ is a primitive group and  $R$ is a $q$-group for some prime 
$q\ne p$. Therefore  $G=R\rtimes M$ for some maximal subgroup $M$ of $G$ 
and  $R=C_{G}(R)=F(G)$ by Theorem 15.2 in \cite[Chapter A]{DH}.

Let   $M_{1}$ be any maximal subgroup of 
$M$. Then $RM_{1}\in F(p)={\cal N}_{p}{\cal H}$.   Since 
  $R=C_{G}(R)$, $O_{q'}(RM_{1})=1$. Hence $O_{q', q}(RM_{1})= 
O_{q}(RM_{1})$  and $O_{p}(RM_{1})=1$.   Therefore   $RM_{1}\in {\cal H}$. 
Let  $H(q)={{\cal G}_{q}}{\cal H}(q)$, where ${\cal H}(q)$ is 
 the  intersection of all formations
containing the set $\{ A/O_{p',\  p}(A) \mid A\in {\cal H}\}$.  Then by 
 Lemma  2.7,  $H(q)=  {\cal   N}_{q}{\cal N}^{r-2}$. Hence  
 $M_{1}/M_{1}\cap RO_{q}(M_{1})   
  \simeq RM_{1}/RO_{q}(M_{1})=RM_{1}/O_{q}(RM_{1})=RM_{1}/O_{q', q}(RM_{1})\in
 {\cal N}_{q}{\cal N}^{r-2}$. 
  Thus  $M_{1}\in  {\cal N}_{q}{\cal 
N}^{r-2}$.  Therefore every maximal subgroup of $M$ belongs to
 $H(q)$.  By induction, ${\cal H}={\cal 
N}^{r-1}$   satisfies   
the boundary condition  in the class of all soluble groups.  Therefore 
 $M\in {\cal H}$, so $G=R\rtimes M\in {\cal F}={\cal N}^{r}$. 
 This contradiction completes the 
proof of the proposition.
 
We use  ${\mathbb P}$  to denote the set of all primes.

{\bf Proposition 4.2.} {\sl Let 
$\{\pi _{i} \mid i\in I   \}$ be a 
 partition of ${\mathbb P}$, and  $ {\cal F}$  the  class of all groups  
$G$   
of the form $G=A_{i_{1}}\times \ldots \times A_{i_{t}}$, where  
$A_{i_{j}}$ is a 
 Hall   ${\pi}_{i_{j}}$-subgroup of $G$,   $i_{1}, \ldots , i_{t}\in I$.  Then 
$ {\cal F}$ is a hereditary saturated formation satisfying the boundary condition. }

{\bf Proof.}   It is clear that the class  $ {\cal F}$  is closed under 
taking subgroups, homomorphic images and direct products. Hence  $ {\cal 
F}$  is a hereditary formation. Moreover, in view of Lemma 2.11 this 
formation  $ {\cal F}$ is saturated.  We show that for any prime $p\in {\pi}_{i}$,
 $F(p)= {\cal G}_{{\pi}_{i}}$. 
 Clearly  
$F(p)\subseteq  {\cal G}_{{\pi}_{i}}$. Suppose that the 
 inverse inclusion is not true and let   $A$ be a group of minimal order 
in $     {\cal G}_{{\pi}_{i}} \diagdown   F(p)$.   Then $A^{F(p)}$ is a unique minimal normal subgroup of 
$A$ and $O_{p}(A)=1$.  Hence   there is a  simple  
 ${\mathbb F}_{p}A$-module $P$    which is     faithful  for $A $ by Lemma 2.6.
Then $G=P\rtimes A \in {\cal G}_{{\pi}_{i}}\subseteq {\cal F}$, so
 $A \simeq G/P=G/O_{p', p}(G) \in F(p)$.  This contradiction shows that  
$F(p)= {\cal G}_{{\pi}_{i}}$. Now let $G$ be any $F(p)$-critical  group.
 Then  $|G|=q$ for some prime $q\not \in {\pi}_{i}$ and so  $G\in {\cal F}$. 
  Hence  ${\cal F}$  satisfies  the boundary condition.

{\bf Proposition 4.3.} {\sl  Let 
$\{\pi _{i} \mid i\in I   \}$ be a 
 partition of ${\mathbb P}$, and  $ {\cal F}$  the  class of all soluble groups  
$G$   
of the form $G=A_{i_{1}}\times \ldots \times A_{i_{t}}$, where  
$A_{i_{j}}$ is a 
 Hall   ${\pi}_{i_{j}}$-subgroup of $G$,   $i_{1}, \ldots , i_{t}\in I$.  Then 
$ {\cal F}$ is a hereditary saturated formation satisfying the boundary condition  
  in the   class of all soluble groups. }

{\bf Proof.}   See the proof of Proposition 4.2.

 {\bf Lattice formations.}    A subgroup $H$ is said to be ${\cal F}$-subnormal 
 in a  group $G$ if either $H=G$ or  there exists a chain of subgroups
 $$H=H_{0} <  H_{1} <     \ldots   <   H_{t}= G $$ such that  
 $H_{i-1}$ is a maximal subgroup of $H_{i}$ and
 $H_{i}/(H_{i-1})_{H_{i}}\in {\cal F}$ for all $i=1,  \ldots , t$
 \cite[p. 236]{Bal-Ez}.

  A 
formation $\cal F$ is said to be  a lattice formation (see \cite[Section 6]{Bal-Ez})
 if the set of all   
${\cal F}$-subnormal subgroups is a sublattice of the lattice of all 
subgroups in every group.

{\bf Proposition  4.4.} {\sl Every  lattice formation   ${\cal F}$ with 
${\cal N}\subseteq {\cal F}\subseteq {\cal S} $   is a  hereditary saturated 
formation  satisfying the boundary condition in the 
class of all soluble groups. }

{\bf Proof.} This follows from Proposition 4.3 and  Corollary 6.3.16 in  \cite{Bal-Ez}. 

{\bf Proposition 4.5.} {\sl   Let ${\cal F}$  be the class of all 
 groups with  $G'\leq F(G)$. Then ${\cal F}$  is a hereditary 
saturated formation satisfying the boundary condition. }

{\bf Proof.} It is clear that ${\cal F}$ is a hereditary formation and ${\cal F}$  is 
saturated by Theorem 4.2 d) in \cite[Chapter III]{Hup}. Moreover, 
 ${\cal F}={\cal N}{\cal A}$, where ${\cal A}$ is the 
 formation of all abelian groups.   Hence by
 Lemma 2.7 (2), $F(p)={\cal G}_{p}{\cal A}$ for
 all primes $p$. Assume that  ${\cal F}$ does not   satisfy  
the boundary condition. Then for some 
prime $p$, the set of all $F(p)$-critical    groups  $A$ with  $A\not \in {\cal F}$
 is non-empty. Let $G $ be a group of minimal order in this set.  Then
  $O_{p}(G)=1=\Phi (G)$ and   $L=G^{\cal F}$ is  a unique minimal normal
 subgroup of $G$ by Lemma 2.10.  Hence  $G$ is a primitive group.

  First we show that 
  $G$ is soluble.  Suppose that this is false. Let $q\ne p$ be any prime
 divisor of $|G|$. Suppose that $G$ is not 
$q$-nilpotent. Then $G$ has a $q$-closed Schmidt subgroup $H=Q\rtimes R$
 \cite[Chapter IV, Satz 5.4]{Hup}, where 
$Q$ is a Sylow $q$-subgroup of $H$, $R$ is a cyclic 
Sylow $r$-subgroup of $H$. Since $G$ is not soluble, 
$H\ne G$. Hence  $H\leq M\in F(p)$ for some  maximal subgroup $M$ of $G$.
 Since $M\in {\cal G}_{p}{\cal A}$,  
  $M'\leq O_{p}(M)$ and hence  $H'\leq Q\cap O_{p}(H)=1$. Therefore  
 $H$ is abelian. This contradiction 
shows that  $G$ is $q$-nilpotent for all primes $q\ne p$, so $G^{\cal N}$
is a Sylow $p$-subgroup of $G$. Hence $G$ is soluble.  
 Therefore
  $L=C_{G}(L)=F(G)$ is a $q$-group for some prime 
$q\ne p$ and   $G=L\rtimes M$ for some maximal subgroup $M$ of $G$ 
 by Theorem 15.2 in \cite[Chapter A]{DH}. 
 Let $M_{1}$ be any maximal subgroup 
of $M$. Then  $LM_{1}\in F(p)$, so $LM_{1}$ is abelian since $L=C_{G}(L)$. 
Hence $M_{1}=1$, so $G'=L$ is nilpotent. Therefore $G\in {\cal F}$.
 This contradiction completes the proof of the 
 result. 

A group $G$ is called a $p$-decomposable if $G=P\times H$, where $P$   is the Sylow 
$p$-subgroup of $G$.

{\bf Corollary 4.6.} {\sl Let  ${\cal F}$  be one of the following  
formations: }

(1) {\sl the class of all nilpotent groups} (Baer \cite{BaerI}); 

(2) {\sl the class of all groups $G$ with $G'\leq F(G)$};

(3) {\sl the class of all $p$-decomposable groups ($p$ is a prime)}.

{\sl Then for each  group $G$,   $Z_{\cal F}(G)=\text{Int}_{\cal 
F}(G)$.}

{\bf Corollary 4.7.} {\sl Let  ${\cal F}$  be one of the following  
formations: }

(1) {\sl the class of all soluble  groups $G$ with $l(G) \leq r \ 
 (r\in {\mathbb N})$} (Sidorov \cite{Sidorov});
 
(2) {\sl any  lattice formation ${\cal F}$ with ${\cal N}\subseteq
 {\cal F}\subseteq {\cal S}$}.

{\sl Then for each  soluble group $G$,  $Z_{\cal F}(G)=\text{Int}_{\cal 
F}(G)$.}

{\bf Some classes of formations not  satisfying the boundary condition.} 
We end this section with some examples of saturated formations which do not satisfy
 the boundary condition.

{ \bf Lemma 4.8 } {\sl 
Let  ${\cal F}$  be any non-empty saturated  formation. 
  Suppose that for some prime $p$ we have $F(p) 
={\cal F}$. Then ${\cal F}$  does not satisfy the boundary condition.}

 {\bf Proof.}  Indeed, in this case  every     ${\cal F}$-critical  
group is also  $F(p)$-critical.

{ \bf Corollary 4.9 } {\sl  Let $p$ be a prime and ${\cal F}$ is one
 of the following formations:}

(1) {\sl  the class of all $p$-soluble   groups;}

(2)    {\sl  the class of all $p$-supersoluble   groups;}

(3)   {\sl  the class of all $p$-nilpotent   groups;}

(4)   {\sl  the class of all soluble   groups.} 

{ Then ${\cal F}$  does not satisfy the boundary condition.}

 {\bf Proof.}  It is clear
 that for any prime $q\ne p$ we have 
 $ {\cal F}={\cal G}_{q}{\cal F}$. Hence   $F(q)={\cal F}$ by Lemma 2.7 
(1).  Now we  use Lemma   4.8.

\section{Further  applications}

Based on the subgroup $ \text {Int}_{\cal F} (G) $ 
   you can achieve the development of many well-known results.
The observations in this section are partial illustrations to this.

{\bf A solubility criterion.}  It is clear that $  \text{Int}_{\cal S}(G)$ is the 
radical $R(G)$ of $G$, that is, the largest soluble normal  subgroup of 
$G$.

{\bf Theorem 5.1.} {\sl  Suppose that a group $G$ has three subgroups
$A_{1}$, $A_{2}$ and $A_{3}$ whose indices $|G:A_{1}|$, $|G:A_2{}|$,
$|G:A_{3}|$ are pairwise  coprime. If $A_{i}\cap A_{j}\leq  R(A_{i})\cap R(A_{j})$
 for all $i\ne j$, then $G$ is soluble.}

{\bf Proof.}   Assume that this theorem  is false and let $G$ be a
counterexample of minimal order.   Fist we shall show that 
 $A_{i}\cap
A_{j}\ne 1$ for all $i\ne j$.   Suppose, for example, that  $A_{1}\cap  
A_{2}=1$. 
Then $A_{1}$ and $A_{2}$ are Hall subgroups of $G$. Hence, for any
prime $p$ dividing $|G:A_{3}|$, $p$ either divides $|G:A_{1}|$  or
divides $|G:A_{2}|$. The contradiction shows that $|G:A_{3}|=1$,
that is, $G=A_3$. Therefore
 $A_{1}$, $A_{2}$ are contained in
$R(G)$. It follows that $G=A_{1}A_{2}=R(G) $, a contradiction.  Therefore
  $A_{i}\cap
A_{j}\ne 1$ for all $i\ne j$.

 Now we 
prove that  $G/N$ is soluble  for any abelian minimal normal
 subgroup $N$ of $G$.  
Let $i\ne j$. Since $N$ is abelian, $N$ is a $p$-group for some prime $p$.
 Hence either  $N\leq A_{i}$ or    $N\leq A_{j}$. In the former case we 
have   
  $$A_{i}/N\cap A_{j}N/N=N(A_{i}\cap A_{j})/N\leq N(R(A_{i})\cap
 R(A_{j}))/N\leq R(A_{i}/N)\cap  R(A_{j}N/N)$$ by Theorem C (a). 
 Therefore the hypothesis holds for $G/N$ 
and so     $G/N$ is soluble   by the choice of $G$.

Finally, we shall prove that $G$ has an abelian minimal normal subgroup.
  Since  $A_{1}\cap A_{2}\ne 1$ and
$A_{1}\cap  A_{2} \leq  R(A_{2})$, for some minimal normal
subgroup $V$ of $A_{2}$ we have  $V\leq  R(A_{2})$. Hence
$V$ is a $p$-group for some prime $p$.  Then either $p$  does not
divide $|G:A_{1}|$ or  $p$ does not divide   $|G:A_{3}|$. Assume
that $p$ does not divide $|G:A_{1}|$. Then for some $b\in A_{2}$, we
have $V\leq A_{1}^{b}$. Hence $V=V^{b^{-1}}\leq A_{1}\cap A_{2}\leq
R(A_{2})$, which implies that
$V^{G}=V^{A_{2}A_{1}}=V^{A_{1}} \leq A_{1}$. It follows that
$E=V^{G}\cap A_{2}\leq A_1\cap A_2\leq R(A_{1})$ and $E$ is
normal in $A_{2}$. Hence $E^{G}=E^{{A_{2}A_{1}}}=E^{A_1}\leq A_1.$
It follows that $E^G=E^{A_1}\leq (R(A_1))^{A_1}=R(A_1)$
 and so $E^G$ is soluble. This shows that $G$
 has an abelian minimal normal subgroup $N$ and we have already proved that $G/N$ is 
soluble, so
 $G$ is soluble contrary to the choice of $G$. This contradiction completes
 the proof of the result.

{\bf Corollary 5.2 (Wielandt \cite{Wiel}).}  {\sl  If $G$ has
three soluble subgroups $A_{1}$, $A_{2}$ and $A_{3}$ whose indices
$|G:A_{1}|$, $|G:A_{2}|$, $|G:A_{3}|$ are pairwise  coprime, then
$G$ is itself soluble. }

{\bf Two characterizations of  supersolubility.}

{\bf Lemma  5.3.}  {\sl Let $N$ be a soluble normal subgroup  of a group 
$G$, $p$  a prime divisor of $|G|$ and $P$ a Sylow $p$-subgroup of $G$.  
Suppose that  $P\not \leq N$ and that every maximal subgroup $M$
of $P$  has a supplement  $T$ in  $G$ such
that $T\cap M\leq \text{Int}_{\cal U }(T)M_{G}$. Then  every maximal subgroup
$V/N$ of $NP/N$ has a supplement  $T/N$
in  $G/N$ such that $$(T/N)\cap (V/N)\leq  \text{Int}_{\cal U }(T/N)(V/N)_{G/N}$$.}

{\bf Proof.} We prove the lemma by induction on $|G|$. Let $V/N$ be
any maximal subgroup of $NP/N$ and $L$ a minimal normal subgroup of
$G$ contained in $N$. Then $L$ is a $q$-group for some prime $q$.
First suppose that $N=L$. If $q\ne p$, then $V=L\rtimes M$, for some
maximal subgroup $M$ of $P$. By hypothesis, 
 there is a subgroup $T$ such that $MT=G$ and $T\cap M\leq
 \text{Int}_{\cal U }(T)M_{G}$.  Then $L\leq T$, $G/L=(V/L) (T/L)$  and 
   $$(V/L)\cap (T/L)=(LM/L)\cap
(T/L)=L(M\cap T)/L    \leq L \text{Int}_{\cal U}(T)M_{G}/L$$  $$
 =(L M_{G}/L)(L\text{Int}_{\cal U}(T)/L)   \leq
  \text{Int}_{\cal U}(T/L) (V/L)_{G/L}     $$  by Theorem C (a).
 If $q=p$, then $V$ is a maximal subgroup of $P$ and so for some 
supplement $T$ of $V$ in $G$ we have   $T\cap V\leq
 \text{Int}_{\cal U }(T)V_{G}$.  Then   $G/L=(V/L)(LT/L)$ and, as above, we
 deduce that  
 $$(V/L)\cap (TL/L)=L(V\cap T)/L\leq \text{Int}_{\cal
U}(T) V_{G}  L/L\leq  \text{Int}_{\cal U}(TL/L)(V/L)_{G/L}$$.

Finally, suppose that $L\ne N$. Obviously, the hypothesis holds for
$(G/L, N/L)$. Hence, by induction, every maximal subgroup
$(V/L)/(N/L)$ of $(PL/L)(N/L)/(N/L)$ 
 has a supplement $(T/L)/(N/L)$ in $(G/L)/(N/L)$
such that 
 $$(T/L)/(N/L)\cap (V/L)/(N/L)\leq 
\text{Int} _{\cal U
}((T/L)/(N/L)) ((V/L)/(N/L))_{(G/L)/ (N/L)}.  $$ Hence from
 the  $G$-isomorphism $G/N\simeq (G/L)/(N/L)$, we
obtain $$(T/N)\cap (V/N)\leq  \text{Int}_{\cal U }(T/N)(V/N)_{G/N}.$$ The
 lemma is proved.

{\bf Theorem 5.4.} {\sl A group $G$ is supersoluble if and only if
every maximal subgroup $V$ of every Sylow subgroup of $G$ 
 has a  supplement $T$ in $G$ such that $V\cap T\leq \text{Int}_{\cal
U}(T)V_{G}$. }
                                       
{\bf Proof.} We need only to prove "if"' part. Suppose that it is false and let G be a
counterexample of minimal order. The proof proceeds via the following steps.

(1) {\sl If $V < P\leq E\leq G$, where  $P$ is a Sylow $p$-subgroup of $G$ and  
  $V$ is a  maximal subgroup  of  $P$, then  
 $V$ has a supplement $T$ in $E$   such that $T\cap V \leq \text{Int}_{\cal
U }(T)V_{E}$}.

Indeed, let $S$ be a supplement of  $V$  in $G$ such that
  $S\cap V \leq \text{Int}_{\cal
U }(S)V_{G}$. Then $T= S\cap E$ is a supplement 
of $V$ in $E$ and  $$V\cap T=V\cap S\cap E \leq 
\text{Int}_{\cal U }(S)V_{G}  \cap E =(\text{Int}_{\cal U }(S)  \cap 
E)V_{G} \leq  \text{Int}_{\cal U }(S\cap E)V_{E} =
 \text{Int}_{\cal U }(T)V_{E}$$ 
  by Theorem C (b).

(2)  {\sl $G/N$ is supersoluble, for every   abelian minimal normal
subgroup $N$ of $G$.}

By Lemma 5.3, the hypothesis is true for $G/N$. Hence $G/N$ is
supersoluble by the choice of $G$.

(3)  {\sl $G$ is soluble. }

In view of (2), it is enough to prove that $G$ has a non-identity
soluble normal subgroup. Suppose that this is false. Then for every maximal
 subgroup $V$ of any Sylow subgroup of $G$ we have $V_{G}=1$.
 Let $p$ be the smallest prime dividing
$|G|$ and $P$ a Sylow $p$-subgroup of $G$. If $|P|=p$, $G$ has a
normal $p$-complement $E$ by \cite[Chapter IV, Theorem 2.8]{Hup}. On the other
hand, by (1), the hypothesis holds for $E$. Hence $E$ is
supersoluble, which implies the solubility of $G$. Hence  $|P| > p$.
If $ V\leq \text{Int}_{\cal U }(G)$ for some maximal subgroup $V$ of $P$, then
$\text{Int}_{\cal U}(G)\ne 1$ and so $G$ has a non-identity soluble normal
subgroup. Therefore  every maximal subgroup $V$
of $P$ has a supplement $T$ in $G$ such that $T\ne G$ and $T\cap
V\leq V_{G}\text{Int}_{\cal U }(T)=\text{Int}_{\cal U }(T)$. We claim that  $T$ is
supersoluble. If $T\cap V=1$, then $|T_p|=p$, for a Sylow
$p$-subgroup $T_{p}$ of $T$. Hence $T$ supersoluble by (1) and the
choice of $G$. Now assume that for some maximal subgroup $V$ of $P$
we have $1\ne T\cap V\leq \text{Int}_{\cal U }(T)$. Since $|P\cap T:V\cap
T|=|V(P\cap T):V|=|P:V|=p$, the order of a Sylow $p$-subgroup of
$T/\text{Int}_{\cal U }(T)$ divides $p$. Hence the hypothesis holds
for $T/\text{Int}_{\cal U }(T)$ by (1) and Lemma 5.3.  
 But since $T\ne G$,
$T/\text{Int}_{\cal U }(T)$ is supersoluble by the choice of $G$. It follows
that $T$ is supersoluble by Theorem C (c). Therefore, our claim holds. This shows
that every maximal subgroup of $P$ has a supersoluble supplement in $G$. By
Lemma 2.12, we see that $G$ has a normal Sylow $q$-subgroup for some prime $q$
 dividing $|G|$.  This contradiction 
completes the proof of (3)

(4)  {\sl $G=N\rtimes M$, where $N=C_{G}(N)=O_{p}(G)$ is a unique  minimal
normal subgroup of $G$ ($p$ is a prime),  $M$ is a supersoluble
maximal subgroup of $G$ with $p$ divides $|M|$ and $|N| > p$.  }

Let $N$ be a minimal normal subgroup of $G$.
 Since  the class of all supersoluble groups is a
 saturated formation, from  (2) and (3) we deduce that  $N$ is
a unique minimal normal subgroup of $G$ and $N\not \leq \Phi (G)$. Hence 
$G$ is a primitive group , so  $N=C_{G}(N)=O_{p}(G)=F(G)$ for some prime 
$p$ by Theorem 15.2 in \cite[Chapter A]{DH}.
  Let $M$ be a maximal subgroup of $G$ such that 
$G=N\rtimes M$. Then $M$ is supersoluble by (2).
 It is also clear that $|N|
> p$. Suppose that $N$ is a Sylow subgroup of $G$ and let $V$ be a maximal
subgroup of $N$. Then $V_{G}=1$, so $V$ has a
supplement $T$ in $G$ such that $VT=G$ and $T\cap V\leq \text{Int}_{\cal
U}(T)$.  But since $T\cap N$ is normal in $G$, the minimality of $N$
implies that either $T=G$ or $T\cap V=1$. In the former case, we
have $1\ne V\leq \text{Int}_{\cal U}(G)$ and so $N\leq \text{Int}_{\cal U}(G)$, which
implies that $G$ is supersoluble by Theorem C (c). In the second case,
 $|T\cap N|=p$, where $T\cap N$
is normal in $G$. Hence $N=N\cap T$ is a group with $|N|=p$. This
contradiction shows that $p | |M|.$ Therefore (4) holds.

(5)  {\sl $\pi (G)=\{p, q\}$, where $p < q$.}

Suppose that $|\pi (G)| > 2$. Let $q\ne p$ be a prime divisor of
$|G|$, $Q$  a Sylow $q$-subgroup of $G$ and  $P$  a Sylow
$p$-subgroup of $G$. Since $G$ is soluble, we may assume that $Q$ and $P$ 
are members of some Sylow system of $G$ and so $E=PQ$ is a proper subgroup of $G$.
 By (1), the hypothesis holds
for $E$. Hence $E$ is supersoluble by the choice of $G$. If $q>p$,
then $Q$ is normal in $E$, which contradicts $C_G(N)=N$. Hence 
 $p>q$ for any prime $q\ne p$ dividing $|G|$.  
 Since $G/N$ is supersoluble, a Sylow $p$-subgroup
$W$ of $G/N$ is normal in $G/N$. Hence $W\leq
O_{p}(G/N).$ By (4), $W\ne 1$. But
$O_{p}(G/N)=O_{p}(G/C_{G}(N))=1$ (see \cite[Appendix, Corollary 6.4]{We}). This
 contradiction shows that $|\pi (G)|=2.$ In the
above proof, we also see that $p>q$ is impossible. Therefore (5)
holds.

{\sl Final contradiction.}  Let $P_{1}$ be a Sylow $p$-subgroup of
$M$, $V$ a maximal subgroup of a Sylow $p$-subgroup $P$ of $G$
containing $P_{1}$ and $M_{q}$ a Sylow $q$-subgroup of $M$.  Then
$N\not \leq V$ and so $V_{G}=1$. By hypothesis, $V$
has a supplement $T$ in $G$ such that $V\cap T\leq
\text{Int}_{\cal U}(T).$ If $T=G$, then $1\ne V\leq \text{Int}_{\cal U}(G)$. Hence
$N\leq \text{Int}_{\cal U}(G)$ since $N$ is the only minimal normal subgroup
of $G$. It follow from (2) that $G$ is supersoluble by Theorem C (c),  a
contradiction. Hence $T\ne G$. In this case, as in the proof of (3),
one can show  that $T$ is supersoluble. Hence a Sylow $q$-subgroup $T_{q}$
of $T$ is normal in $T$ by (5). But $T_{q}$ is a Sylow subgroup of $G$.
Hence $T_{q}=(M_{q})^{x}$ for some $x\in G$. Since $q > p$ and $M$
is supersoluble,  $M= N_{G}(M_{q})$. Hence $T\leq 
N_G(T_q)=N_G(M_q^x)=(N_G(M_q))^x=M^{x}$. But then $G=VT=VM^{x}=VM$
by Lemma 2.13. It follows that $|G|=|V||M|/|V\cap M|\leq
|V||M|/|P_{1}| <  |M||N|=|G|$. This contradiction completes the
proof of the result.

 Note that if $H$ is a group of $G$ and  $H$ either is normal in $G$, has a 
complement in $G$, or has a supplement $E$ in $G$ with $E\in {\cal F}$, then  
$H$  has a  supplement $T$ in $G$ such that $V\cap  
 T\leq \text{Int}_{\cal U}(T)V_{G}
$. Hence from Theorem 5.4 we get  the following

{\bf Corollary 5.5 (Srinivasan \cite{5})}. {\sl
If the maximal subgroups of the
Sylow subgroups of $G$ are normal in $G$, then $G$ is supersoluble.

}

{\bf Corollary 5.6 (Ballester-Bolinches and  Guo \cite{B-G}).}
{\sl A group $G$ is supersoluble if  every maximal
subgroup of every Sylow subgroup of $G$  has a
complement  in $G$.}

{\bf Corollary 5.7 (Guo, Shum and  Skiba \cite{GSS}).}
{\sl A group $G$ is supersoluble if and only if every maximal
subgroup of every Sylow subgroup of $G$ has a
supersoluble supplement in $G$.}

In view of Theorem C (h) for every group $G$ we have  
 $Z_{\cal F}(G)\leq \text{Int}_{\cal F}(G)$. Hence from  Theorem 5.4 we also 
get 

{\bf Corollary 5.8 (Guo and  Skiba \cite{GSJGT}).} {\sl 
 A group $G$ is supersoluble if and only if every maximal subgroup $V$ of
every Sylow subgroup of $G$ either is normal or has a supplement $T$ in $G$
 such that   $V \cap T \leq  Z_{\cal U}(T)$.  }

It is well known that   if every minimal  subgroup of a group   $G$  is
normal  in $G$, then the commutator subgroup $G'$  of $G$ is 2-closed
  (Gasch\"utz \cite[IV, Theorem 5.7]{Hup}). On the other hand, if $G$ 
is a group  of odd order and  every minimal subgroup
 of $G$ is  normal in $G$, then $G$ is
supersoluble (Buckley \cite{3}).  The following theorem covers both these 
observations.

{\bf Theorem 5.9.} {\sl A group  $G$  is $2'$-supersoluble if and only if 
 every minimal subgroup $L$ of  $G$ of odd order is contained in the
 intersection of
 all maximal $2'$-supersoluble subgroups of $G$. }

{\bf Proof.}    Let    ${\cal F}$  be the class of all
 $2'$-supersoluble groups  and $I=\text{Int}_{\cal F}(G)$
  the intersection of all maximal $2'$-supersoluble subgroups of 
$G$. It is well known that  the class   ${\cal F}$  is a hereditary saturated 
formation (see \cite[Chapter VI, Satz 8.6]{Hup}).   Assume that every 
minimal subgroup $L$ of  $G$ of odd order is contained in $I$.  We shall prove 
that $G$ is  $2'$-supersoluble.  Assume that this   
 is false and let $G$ be a
counterexample of  minimal order .

The hypothesis holds for
every subgroup of $G$ by Theorem C (b). Hence every maximal subgroup  of $G$ is 
$2'$-supersoluble  by the choice of $G$. Therefore every maximal subgroup of $G$
is soluble.

  First we show that $G$ is soluble. Assume that
this  is false. Then $G=G'$, and if  $F=F(G)$, then $F=\Phi (G) $,  
 $G/F$ is a simple
non-abelian group and every proper normal subgroup of $G$ is contained in
$F$. Hence $I=F$. It is clear that every maximal subgroup of $G/F$
 is soluble and hence by \cite{Tomph},  $G/F$ is isomorphic to one of the
 following groups:
 $PSL_{2}(p)$ (where  $p >3 $ is a prime
 such that    $p^{2}+1 \equiv 0(5) $), 
 $PSL_{2}(3^{p})$ (where $p$ is an odd prime), $PSL_{2}(2^{p})$ (where 
 $p$ is a prime),  $PSL_{3}(3)$, a  Suzuki group
$Sz(2^{p})$ (where $p$ is an odd prime).

   Let $r$ be the  largest prime dividing $|G/F|$ and $G_{r}$ a Sylow 
$r$-subgroup of $G$. Then $r >3$ by Burnside's $p^{a}q^{b}$-theorem. 
Let $p$ be any   odd prime dividing $|G/F|$ and  $C_{p}$ a
 subgroup of $G$ of order   $p$.     
 Then  $C_{p}\leq I=F$. Suppose that $p < r$ and let $P$ be a Sylow 
$p$-subgroup of $F$. We show that   $E=P{G_{r}}^{x}$ is 
nilpotent  for all $x\in G$.   Suppose that this is false and
 let $H$ be  a Schmidt subgroup 
of $E$, that is,  an  ${\cal N}$-critical  group. 
  Since $G$ is not soluble, $E\ne G$ and 
hence $H$ is supersoluble. Therefore ${G_{r}}^{x}$ is normal
 in $H=P\rtimes {G_{r}}^{x}$  since    $p < r$, so $H$ is nilpotent. 
 This contradiction shows that   $P{G_{r}}^{x}$ is 
 nilpotent. Hence $\langle (G_{r})^{G}\rangle =G\leq
C_{G}(P)$.  Thus $P\leq Z(G)$ and  $P\leq \Phi (G)$ since   $F=\Phi (G)$.   
Let $V$ be a Hall $p'$-subgroup of $F$. Then   $PV/V\leq Z(G/V)$ and  
$PV/V\leq \Phi (G/V)$. Hence 
 $p$ divides $|M(G/F)|$, where $M(G/F)$
 is the Schur multiplicator of $G/F$. 
 Since $p > 2$, it follows that 
$p=3$,  $\pi (|M(G/F)|)\subseteq \{2, 3\}$   and $5$ 
divides $|G/F|$ (see \cite[Chapter 4]{GorII}). 
 Let   
  $G_{3}$ be a 
Sylow $3$-subgroup of $G$ and $R$  the Sylow $5$-subgroup of $F(G)$.  
Since   $V=RG_{3}$ is soluble, $V\ne G$ and so  $V$ is supersoluble. 
 Hence for any chief factor $H/K$ of $V$ below   
$R$ we deduce  that $|V/C_{V}(H/K)|$ divides $4$. Therefore $C_{V}(H/K)=V$, so 
$R\leq Z_{\infty}(V)$ and hence $V$ is nilpotent.
Thus  $R\leq Z(G)$, which implies that $5$ divides  $|M(G/F)|$,  a 
contradiction. Therefore $G$ is soluble. But $G$ is an ${\cal F}$-critical  group.
  Hence by Lemma 2.8, $G^{\cal F}$ is a $p$-group  for some odd prime $p$
and 
 $\psi _{e} (G^{\cal F})\leq I$. Thus 
$G\in {\cal F}$  by Theorem C (c)(g). This contradiction completes 
the proof of the result.

{\bf A nilpotency criterion.}  
In the following theorem,  $c(G)$ denotes the nilpotent class of the 
nilpotent group $G$.

{\bf Theorem 5.10.} {\sl   Suppose that $G$ has three subgroups
$A_{1}$, $A_{2}$ and $A_{3}$ whose indices $|G:A_{1}|$, $|G:A_2{}|$,
$|G:A_{3}|$ are pairwise coprime. Suppose that
 $A_{i}\cap A_{j}\leq  Z_{n}(\text{Int}_{\cal N }(A_{i}))\cap
Z_{n}(\text{Int}_{\cal N}(A_{j}))$ for all $i\ne j$. Then  $G$ is nilpotent
and  $c(G) \leq n$.}

{\bf Proof.}  Let
 $p$ be any prime dividing $|G|$.
 By
hypothesis, there exists $i\neq j$  such that $p\nmid |G:A_i|$ and
$p\nmid |G:A_j|.$ Hence $p\nmid |G:A_i\cap A_j|$ and so $G$ has a Sylow
$p$-subgroup $P$ such that
 $P\leq A_{i}\cap A_{j}\leq
Z_{n}({\text{Int}} _{\cal N}(A_{i}))\cap Z_{n}({\text{Int}} _{\cal N}(A_{j}))$. 
 Since 
    $ \text{Int}_{\cal N}(A_{i})$ is nilpotent, 
  $P$ is a
characteristic subgroup of $\text{Int} _{\cal N}(A_{i})$. On the other hand, 
 $ \text{Int}_{\cal N}(A_{i})$ is characteristic in  $A_{i}$. Hence   $P$ is
normal in $A_{i}$.  Similarly, we have $A_{j}\leq N_{A_{j}}(P)$.
Therefore  $G=A_{i}A_{j}\leq N_{G}(P)$. Thus $G$ is nilpotent and 
$c(P)\leq n$ for all Sylow subgroups $P$ of $G$, which implies $c(G)\leq n$.

The theorem is proved.

{\bf Corollary 5.11 (Kegel \cite{Kegel})}. {\sl  If $G$ has
three nilpotent subgroups $A_{1}$, $A_{2}$ and $A_{3}$ whose indices
$|G:A_{1}|$, $|G:A_{2}|$, $|G:A_{3}|$ are pairwise  coprime, then
$G$ is itself  nilpotent. }

{\bf Corollary 5.12
 (Doerk \cite{Do})}. {\sl  If  $G$ has three
abelian subgroups $A_{1}$, $A_{2}$ and $A_{3}$ whose indices
$|G:A_{1}|$, $|G:A_{2}|$, $|G:A_{3}|$ are pairwise  coprime, then
$G$ is itself abelian. }

{\bf A question of Agrawal.}   
 Recall that a subgroup
  $H$ of a group $G$ is said to be $S$-quasinormal in $G$
 if $HP=PH$ for all Sylow subgroups 
$P$ of $G$. The 
 hyper-generalized-center $genz^{*}(G)$  of  $G$ coincides with the largest term
 of the chain   of subgroups

 $$1=Q_{0}\leq Q_{1}\leq  \ldots \leq  Q_{t}\leq \cdots $$
where $Q_{i}(G)/Q_{i-1}(G)$ is the subgroup of $G/Q_{i-1}(G)$ generated by 
the set of all cyclic $S$-quasinormal subgroups of  $G/Q_{i-1}(G)$ (see
 \cite[page 22]{We}). In the paper  \cite{Agr}, Agrawal proved 
 that  $genz^{*}(G)$ is contained in every 
maximal supersoluble subgroup of the group $G$ and  posed the following 
question:  {\sl  Does there exist a group  $G$ with
 $genz^{*}(G) \ne  \text{Int}_{\cal
U}(G)$? }
 (see \cite[page 19]{Agr} or \cite[page 22]{We})

The following  example gives  a positive answer to this question and shows that
 there are soluble groups $G$ with $\text{Int}_{\cal  U}(G)\ne Z_{\cal U}(G)$.

{\bf Example 5.13.}  Let $ C_{p}$ be a group of prime order $p$ with
 $|\pi (\text{Aut}(C_{p}))| > 1$. Let  $R$ and $L$  be Hall subgroups of
 $\text{Aut}(C_{p})$
such that   $\text{Aut}(C_{p})=R\times L$ and for any 
 $r \in \pi    (R)$ and  $q\in \pi   (L) $ we have $ r <  q$.
Let  $ G =(C_{p}\rtimes R) \wr L=K \rtimes  L $ be the regular 
wreath product of $C_{p}\rtimes R$  with  $L$, where  $K$ is the base 
group of $G$.
 Let $P={C_{p}}^{\natural}$ (we use here the terminology in \cite[Chapter A]{DH}).
Then by Proposition 18.5 in \cite[Chapter A]{DH},  $G$ is a primitive 
group and    $P$ is a unique minimal  normal subgroup  of $G$. Hence 
$P=F(G)=C_{G}(P)$. Moreover, by Lemma 18.2 in   \cite[Chapter A]{DH}, 
$G=P\rtimes  M$, where $M\simeq U= R \wr L=D\rtimes L$, where $D$ is 
 the base group of 
$U$. It is clear that  $D$ is a Hall abelian subgroup of $U$ 
and $L$ is a cyclic  subgroup of $U$ such that  
 for any        $r \in \pi    (D)$ and  $q\in \pi   (L) $ we have $ r <  
q$. Moreover, since $|\text{Aut}(C{p})|=p-1$, $D$ and $L$ are groups of 
exponent dividing $p-1$.  
First we show that every supersoluble subgroup $W$ of $U$ is 
nilpotent. Suppose that this is false and let $H$ be a Schmidt subgroup of 
$W$. Then 
 $1 < D\cap H < H$, where  $D\cap H$ is a Hall normal subgroup of 
$H$.    By \cite[Chapter IV, Satz 5.4]{Hup}, there are primes $r$ and $q$
such that  $H=H_{r}\rtimes H_{q}$, where 
$H_{r}$ is a Sylow $r$-subgroup of $H$, $H_{q}$ is a cyclic Sylow $q$-subgroup of
 $H$. 
 Hence $D\cap H=H_{r}$. Since  $H\leq W$, 
$H$ is supersoluble and hence $r > q$. But  $Q\simeq H/D\cap H \simeq HD/D$ is 
isomorphic with some subgroup of $L$, so $r < q$. This contradiction shows that $W$
 is nilpotent.

Now we shall show that  $P\leq \text{Int}_{\cal U}(G)$. Let $V$ be any 
supersoluble subgroup of $G$ and $ W$   a Hall $p'$-subgroup of $ V$.
 Then $ PV = P W$. It is clear that $M$ is a Hall $ p'$-subgroup of $ G$, so  
for some $x\in G$ we have $W\leq  M^{x}\simeq  U^{x}$.  Hence $ W$ is nilpotent
 since $W$ is a subgroup of the supersoluble group $V$.   
 It is clear that the Sylow subgroups of $W$ are abelian, hence
 $W$ is an
 abelian group of exponent dividing $p-1$. Hence $ PV $ is supersoluble by 
\cite[Chapter 1, Theorem 1.9]{We}.  Therefore $P\leq \text{Int}_{\cal 
U}(G)$.

Finally, we show that $ genz^{*}(G)=1$. Indeed, suppose that  $ 
genz^{*}(G)\ne 1$.   Then $G$ has a non-identity cyclic  $S$-quasinormal 
subgroup, say $V$. The subgroup $V$ is subnormal in $G$ by 
 \cite[Chapter  1, Corollary 6.3]{We}. Therefore   
 $ V \leq  F(G)=  P  $ by \cite[Chapter A, Theorem 8.8]{DH}.  Moreover, if 
$Q$ is a Sylow $q$-subgroup of $G$, where $q\ne p$, then $V$ is subnormal 
in $VQ$ and so $Q\leq N_{G}(V)$. Hence $V$  is normal
 in $G$ and therefore $ V =  P$ is cyclic. But then $|P|=p=|{C_{p}}^{\natural}|$,
 a contradiction.     
Thus  $  genz^{*}(G)=1= Z_{\cal U}(G)$.

{\bf Acknowledgment.}    The author is very  grateful for the helpful suggestions
 of the referee. The author is  also indebted to  Professor  L.A. Shemetkov  
for his useful suggestions and  comments. 

 particular

\end{document}